\documentclass[preprint,12pt]{amsart}
\usepackage{graphicx}
\usepackage{amssymb}
 \usepackage{amsthm}
 \usepackage{amsmath}
 \usepackage{amssymb}
 \usepackage{float}
\usepackage{caption}
\usepackage{subcaption}
\usepackage[numbers]{natbib}
\usepackage[colorlinks=false]{hyperref}
\usepackage{color}
\usepackage{natbib}

\title{A Modified Quadratic Lorenz Attractor}

\author{Bu\u{g}\c{c}e Emina\u{g}a}
 \address{Department of Computer Engineering, Girne American University, Kyrenia, North Cyprus, via Mersin 10, Turkey} 
 \email{bugceeminaga@gau.edu.tr}
 \author{Hatice Akt\"ore }
 \address{Department of Mathematics, Eastern Mediterranean University, Famagusta, North Cyprus, via Mersin 10, Turkey}
 \email{hatice.aktore@emu.edu.tr}
 \author{Mustafa Riza}
 \address{Department of Physics, Eastern Mediterranean University, Famagusta, North Cyprus, via Mersin 10, Turkey}
 \email{mustafa.riza@emu.edu.tr}


\begin{document}
 \maketitle

 \begin{abstract}
 This study introduces a modified quadratic Lorenz attractor. The properties of this new chaotic system are analysed and discussed in detail, by determining the equilibria points, the eigenvalues of the Jacobian,  and the Lyapunov exponents. The numerical simulations, the time series analysis, and the projections to the $xy$-plane, $xz$-plane, and $yz$-plane  are conducted to highlight the chaotic behaviour. The multiplicative form of the new system is also presented and the simulations are conducted using multiplicative Runge-Kutta methods.
 \end{abstract}

\section{Introduction}

Dynamical systems are  mathematical models describing the evolution of systems in terms of equation of motion with sensitive initial values. The application areas of dynamical systems can be found e.g. in many disciplines in physics \cite{roessler76,LR,arnold}, population dynamics in biology \cite{Bio}, and chemical kinetics in chemistry \cite{chemical1,chemical2}. Dynamic systems theory has also found its way into subjects outside the fields of Mathematics and natural sciences, exemplarily we would like to mention mathematical economy and finance \cite{MCRO,EC,peter}. Poincar{\'e} made a kick start to the subject of dynamical systems in his pioneering work  in 1890 \cite{P1}. Later, in the 1920's, \citet{fatou17} and \citet{julia18} introduced the dynamics of complex analytic maps. On the  one hand,  studies like \citet{BGD1}, \citet{K1,K2,K3}, Cartwright and Littlewood \cite{cartwright} find their origin  in physics, like e.g. the three body problem in astronomy, whereas on the other hand  Stephan Smale provided a purely mathematics motivated approach  \cite{ss1,ss2}. E. N. Lorenz \cite{LR} observed that very simple differential equations  become chaotic under certain circumstances. The system proposed by Lorenz shows a very complex dynamical behaviour and displays the well-known two-scroll butterfly-shape. The dynamic equations of the Lorenz system are given as
\begin{eqnarray}
\frac{dx}{dt}&=&\sigma(y-x), \label{eq:lorx}\\
\frac{dy}{dt}&=&x(\rho-z)-y, \label{eq:lory}\\
\frac{dz}{dt}&=&xy-\beta z, \label{eq:lorz}
\end{eqnarray}
where the parameters $\sigma$, $\rho$, and $\beta$ are assumed to be positive. Lorenz used exemplarily the values $\sigma = 10$, $\beta = 8/3$ and $\rho = 28$ to demonstrate the systems chaotic behaviour. The study on "The equation for continuous chaos" by R\"ossler \cite{roessler76} can be seen as another landmark in the discussion of 3D dynamic systems. The R\"ossler system and the Burke Shaw  system \cite{shaw81} have the property of two unstable saddle foci in common. Other chaotic systems exhibiting a similar simple structure as the Lorenz system, without being topologically equivalent, are proposed by Chen \cite{chen99} and L\"u  and Chen \cite{Lu}. \cite{Lu} discusses the transition between Lorenz and Chen attractors. Furthermore, Yang and Chen introduced another chaotic system with three fixed points: one saddle and two stable fixed points. Yang et al. \cite{Y13} and Pehlivan et al. \cite{I4} introduced and analysed  chaotic systems  similar to the Lorenz, Chen, and Yang-Chen systems, with two different fixed points, i.e. two stable node-foci. Chaotic systems have found their ways also into many applications in engineering, such as electronic circuits \cite{I5}-\cite{I20}. Modern pacemakers actually rely nowadays on these chaotic circuits. 
    
    Dynamical systems have also been discussed in the framework of various Non-Newtonian Calcului, like fractional calculus, geometric multiplicative calculus, and bigeometric multiplicative calculus. Jun Guo Lu transformed the L\"{u} system into fractional calculus  \cite{fLu}, investigating  chaotic behaviour of  fractional-order of the L\"{u} system numerically. Aniszewska applied  multiplicative  calculus to the R\"{o}ssler system and showed chaotic behaviour of multiplicative R\"{o}ssler\cite{A}. As an application of the bigeometric multiplicative Runge-Kutta method, the bigeometric multiplicative R\"ossler System was solved in \cite{MB}.

    
 This study introduces a new chaotic attractor, found by modification of the Lorenz system by a quadratic term.  Detailed numerical and theoretical analysis reveals that the proposed system shows chaotic behaviour and the property of a two-scroll attractor like the Lorenz attractor. 
  
Section \ref{sec:1} introduces the modified quadratic Lorenz attractor, discusses and analyses its properties by determining the equilibria points, and Lyapunov exponents theoretically as well as numerically.  In Section 3, the modified quadratic Lorenz attractor is translated into  geometric and bigeometric calculus, and the solutions of the the system are obtained using the corresponding multiplicative Runge-Kutta methods \cite{RA13} and \cite{MB}. Finally, this paper closes with the summary of the obtained results.

\section{Design of a new Chaotic System}
\label{sec:1}
This paper presents a new Chaotic system derived from the Lorenz system. The system is generated by the following simple three-dimensional system:

\begin{eqnarray}
\frac{dx}{dt}&=&\sigma(y z-x) \label{eq:nlorx}\\
\frac{dy}{dt}&=&\rho x-x z  \label{eq:nlory}\\
\frac{dz}{dt}&=& (x y)^2-\beta z  \label{eq:nlorz} 
\end{eqnarray}

where $x$, $y$, and $z$ are variables and $\sigma$, $\rho$, and $\beta$ are real parameters. In the new proposed chaotic system, all  the equations have some differences compared to the original Lorenz system. In order to see the differences between the two systems we can compare the equations  \eqref{eq:lorx}-\eqref{eq:lorz} and  \eqref{eq:nlorx}-\eqref{eq:nlorz} one by one.  Evidently, in the Lorenz system equation \eqref{eq:lorx} is linear, whereas  equation \eqref{eq:nlorx} is non-linear.  Furthermore, equations \eqref{eq:lory} and \eqref{eq:nlory} are both non-linear, where the $y$-dependence of  equation \eqref{eq:nlory} is cancelled. The most significant difference can be observed comparing equations \eqref{eq:lorz} and \eqref{eq:nlorz}. In \eqref{eq:nlorz} the term $xy$ is squared compared to equation \eqref{eq:lorz}.

\subsection{System Description}

\begin{figure}[H]
    \centering
    \begin{subfigure}[b]{0.4\textwidth}
        \centering
        \includegraphics[width=\textwidth]{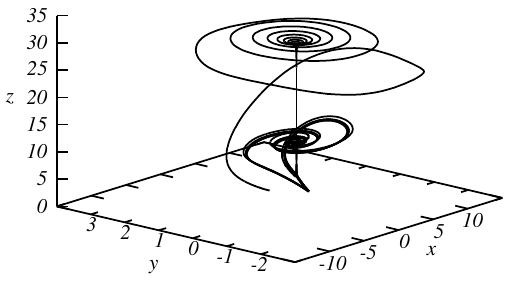}
        \caption{$\beta=0.1$}
        \label{fig:y equals x}
    \end{subfigure}
    \hfill
    \begin{subfigure}[b]{0.4\textwidth}
        \centering
        \includegraphics[width=\textwidth]{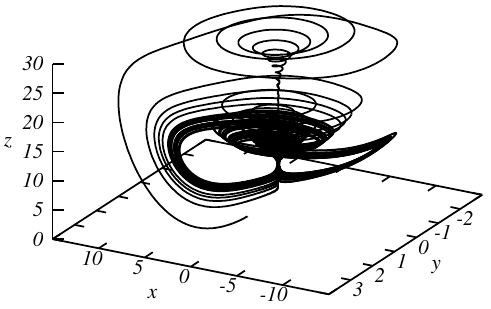}
        \caption{$\beta=0.5$}
        \label{fig:three sin x}
    \end{subfigure}
    \hfill
    \begin{subfigure}[b]{0.4\textwidth}
        \centering
        \includegraphics[width=\textwidth]{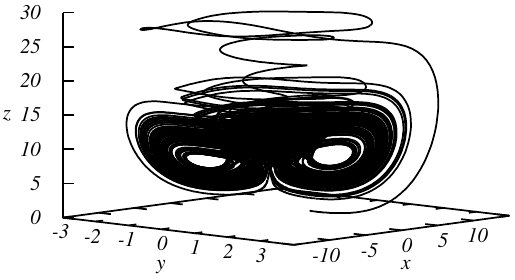}
        \caption{$\beta=2$}
        \label{fig:five over x}
    \end{subfigure}
 \hfill
    \begin{subfigure}[b]{0.4\textwidth}
        \centering
        \includegraphics[width=\textwidth]{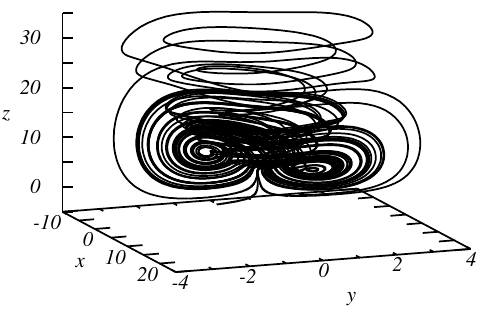}
        \caption{$\beta=4$}
        \label{fig:five over x}
    \end{subfigure}
 \hfill
    \begin{subfigure}[b]{0.4\textwidth}
        \centering
        \includegraphics[width=\textwidth]{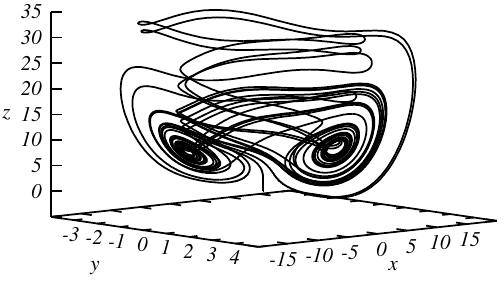}
        \caption{$\beta=10$}
        \label{fig:five over x}
    \end{subfigure}
    \caption{Simulation of the system when  $\sigma=12$, $\rho=8$ and various parameter values }
    \label{fig:three graphs}
\end{figure}

The initial values and the parameters of the system are chosen as $(1,1,1)$ and $\sigma=12$, $\rho=8$ with varying $\beta$. Then it has been observed that the behavior of the new chaotic system changes depending on the different values of $\beta$ as shown in Figure \ref{fig:three graphs}.

As a result of the calculations, as it can be seen in Figure \ref{fig:three graphs}, it has been observed that the new system is  chaotic for the  parameters 
\begin{equation}
\sigma=12,  \rho=8, \text{ and } \beta=4.  \label{param}
\end{equation}
The analysis of the new chaotic system will be done according to those parameters. 

\subsection{System Analysis}

The first step to analyze a chaotic system is to find the equilibrium points. In order to determine the equilibrium points of the proposed system \eqref{eq:nlorx}-\eqref{eq:nlorz}, we need to solve the system

\begin{eqnarray}
\sigma(yz-x)&=&0 \nonumber\\
\rho x-xz&=&0 \label{eq:equi} \\
(xy)^2-\beta z&=&0\nonumber
\end{eqnarray}
 
Thus, the solution of the system \eqref{eq:equi} with respect to $x$, $y$, $z$ give the equilibria points as:

\begin{eqnarray*}
\it{O}&=&(0,0,0) \\
\it E^+&=&\left(\sqrt[4]{\beta\,\rho^3},\sqrt[4]{\frac{\beta}{\rho}},\rho\right) \\
\it E^-&=&\left(-\sqrt[4]{\beta\,\rho^3},-\sqrt[4]{\frac{\beta}{\rho}},\rho\right)
\end{eqnarray*}

For the parameters  chosen in \eqref{param}, the numerical values of the equilibria points are 

\begin{eqnarray}
\it{O}&=&(0,0,0), \label{eq:eqi1}\\
\it E^+&=&(6.73,0.84,8), \label{eq:eqi2}\\
\it E^-&=&(-6.73,-0.84,8).\label{eq:eqi3}
\end{eqnarray} 

In order to decide on the stability of the new proposed system, the eigenvalues of the Jacobian matrix must be analyzed. The Jacobian matrix for this system  \eqref{eq:nlorx}-\eqref{eq:nlorz} can be easily obtained as
\begin{equation}
\it{J} = 
\begin{bmatrix}
      -\sigma &  \sigma z &  \sigma y \\
       \rho -z & 0           & -x \\
       2xy^2           & 2x^2y & -\beta
     \end{bmatrix}.\label{eq:jacobian}
\end{equation}
 
The expressions for the eigenvalues of the Jacobian matrix  \eqref{eq:jacobian} are very long and complicated. As we are only interested in the numerical values of the eigenvalues at the equilibria points \eqref{eq:eqi1}-\eqref{eq:eqi3} for the given parameters  \eqref{param}, the eigenvalues are stated in the table below:
\begin{table}[h]
\begin{center}
\begin{tabular}{|l|c|c|c|}
\hline
Equilibrium & $\displaystyle \lambda_1$& $\displaystyle \lambda_2$& $\displaystyle \lambda_3$ \\
Point  &&&\\
\hline
$\displaystyle \it{O}$ & $\displaystyle -12$ & $\displaystyle -4$ &$\displaystyle 0$ \\ \hline
$\displaystyle \it E^+$ & $\displaystyle 2.65+23.87 i$ & $\displaystyle 2.65-23.87 i$ &$\displaystyle-21.3$ \\ \hline
$\displaystyle \it E^-$ & $\displaystyle 2.65+23.87 i$ & $\displaystyle 2.65-23.87 i$ &$\displaystyle-21.3$ \\ \hline
\end{tabular}
\end{center}
\caption{Eigenvalues of the Jacobian at the equilibrium points}
\end{table}

As the eigenvalues $\lambda_1$ and $\lambda_2$ for the equilibrium point $\it O$ are both negative, the system is unstable at this equilibrium point.  
 The eigenvalues corresponding to the equilibrium point $\it E^-$ will be the same with the eigenvalues of $\it E^+$, because of the quadratic nature of the system. Since $\lambda_3$ is a negative real number and $\lambda_1$ and $\lambda_2$ are two complex conjugate eigenvalues with positive real parts, equilibrium points  $\it E^+$ and $\it E^-$ are unstable according to \cite{St}. 

\subsection{Symmetry and Dissipativity}

The System \eqref{eq:nlorx}-\eqref{eq:nlorz} has a natural symmetry and is invariant under the coordinate transformation $(x,y,z)\rightarrow(-x,-y,z)$ which persists for all values of the system parameters. So, system \eqref{eq:nlorx}-\eqref{eq:nlorz} has rotation symmetry about the $z-axis$.

Let, $f_1=\frac{dx}{dt}$ , $f_2=\frac{dy}{dt}$ and $f_3=\frac{dz}{dt}$ in the system \eqref{eq:nlorx}-\eqref{eq:nlorz}. Then we get for the vector field

\[
(\dot x, \dot y, \dot z)^T= \left(f_1, f_2, f_3\right)^T
\]

Consequently  the divergence of the vector field $\mathbf V$ yields to:

\begin{equation}
{\mathbf {\nabla}} \cdot (\dot x, \dot y, \dot z)^T= \frac{\partial f_1}{\partial x}+ \frac{\partial f_2}{\partial y}+\frac{\partial f_3}{\partial z}=-(\sigma+\beta)=f.
\label{eq:diss}
\end{equation}

Note that $ f = -(\sigma+\beta) =-16$  is a negative value, so the system is a dissipative system and an exponential rate is:

\begin{equation}
\frac{dV}{dt}= f V  \Longrightarrow  V(t) = V_0 e^{f t}=V_0e^{-16 t}
\label{eq:exprate}
\end{equation}

From \eqref{eq:exprate}, it can be seen that a volume element $V_0$ is contracted by the flow into a volume element $V_0e^{-16t}$ at the  time $t$ .

\subsection{Lyapunov Exponent and Fractional Dimension}

The Lyapunov exponents generally refer to the average exponential rates of divergence or convergence of nearby trajectories in the phase space. The important part is that if there is at least one positive Lyapunov exponent, the system can be defined to be chaotic. According to the detailed numerical and theoretical analysis, the Lyapunov exponents are found to be $l_1 =5.4162$,  $l_2 = 2.1912$, and $l_3 = -19.2269$. Therefore, the Lyapunov dimension of this system is:

\begin{equation}
D_L=\displaystyle{j+\displaystyle\frac{\displaystyle\sum_{i=1}^j l_i}{\left|l_{j+1}\right|}=2+\frac{l_1+l_2}{\left|l_3\right|}}=2.3957
\label{eq:dim}
\end{equation}
 This result is consistent with the findings in \cite{lyacompare}, i.e. that the Lyapunov dimension is in the range 2-3.
Equation \eqref{eq:dim} means that the system \eqref{eq:nlorx}-\eqref{eq:nlorz} is a dissipative system, and the Lyapunov dimensions of the system are fractional. Having a strange attractor and positive Lyapunov exponent, it is obvious that the system is a 3D chaotic system.

\begin{figure} [H]
\centerline{\includegraphics[width=10cm]{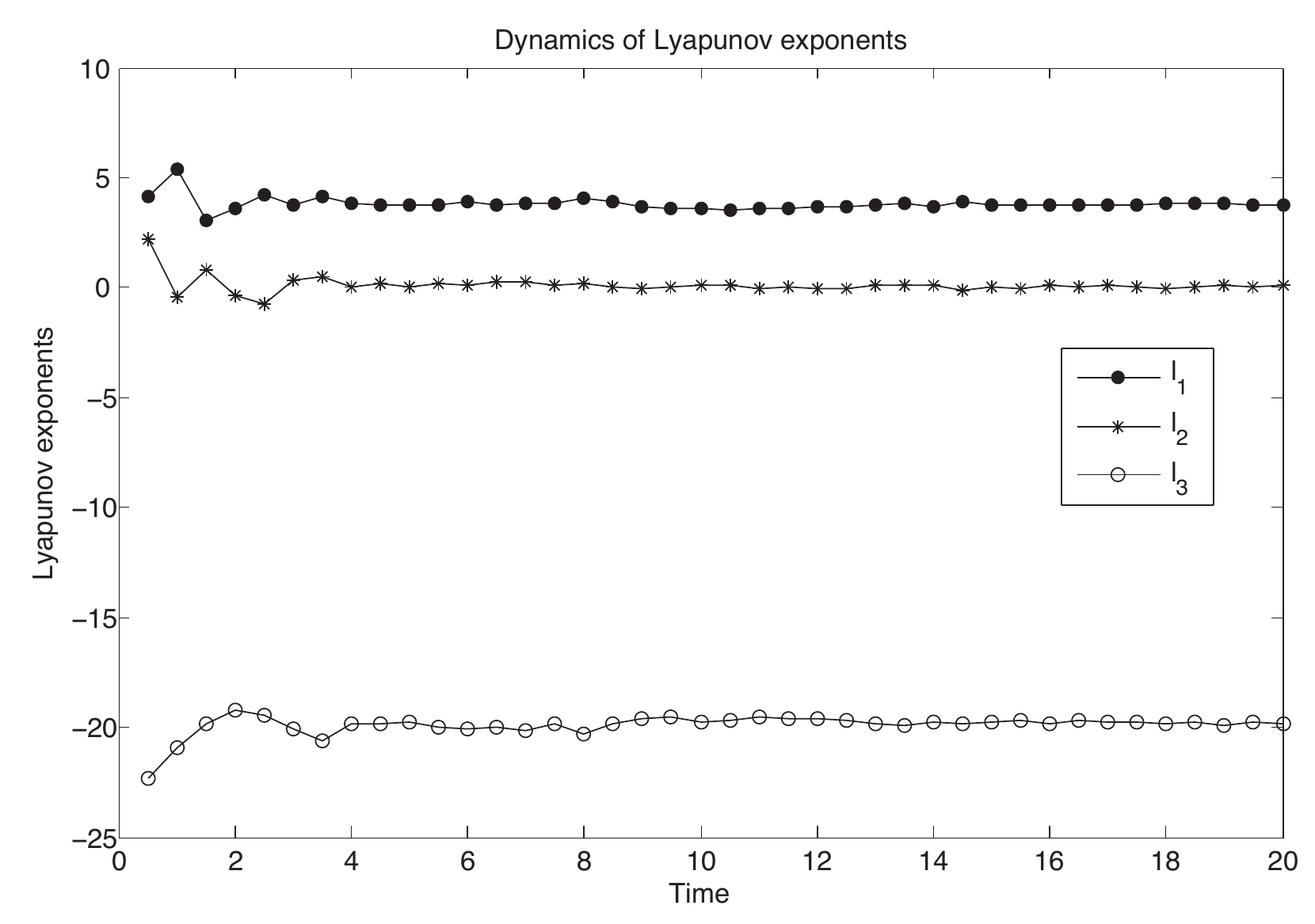}}
\caption{Plot of Lyapunov exponents}
\end{figure}

\subsection{Numerical Simulations}

\begin{figure}[H]
    \centering
    \begin{subfigure}[b]{0.4\textwidth}
        \centering
        \includegraphics[width=\textwidth]{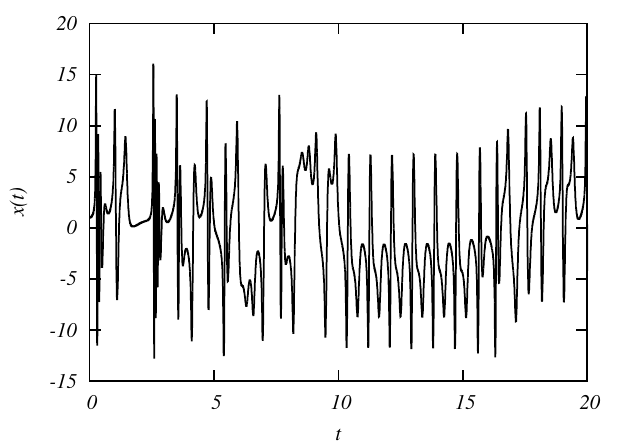}
        \caption{}
        \label{fig:x(t)}
    \end{subfigure}
    \hfill
    \begin{subfigure}[b]{0.4\textwidth}
        \centering
        \includegraphics[width=\textwidth]{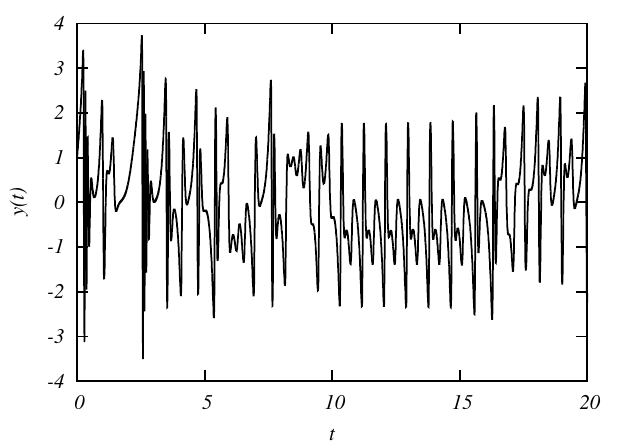}
        \caption{}
        \label{fig:y(t)}
    \end{subfigure}
    \hfill
    \begin{subfigure}[b]{0.4\textwidth}
        \centering
        \includegraphics[width=\textwidth]{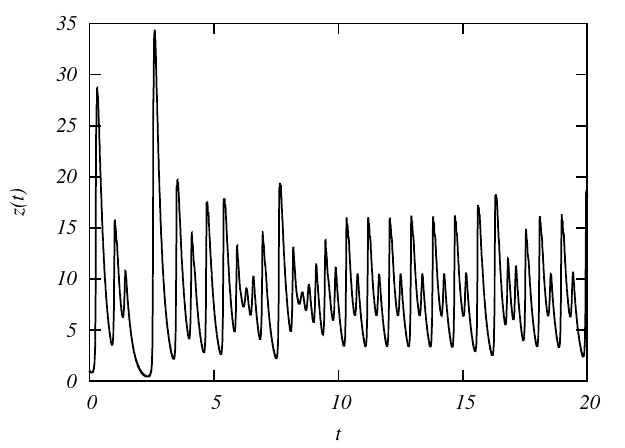}
        \caption{}
        \label{fig:z(t)}
    \end{subfigure}
    \caption{Waveforms of $x(t)$, $y(t)$, $z(t)$ respectively}
    \label{fig:wf}
\end{figure}

Time series analysis of the system \eqref{eq:nlorx}-\eqref{eq:nlorz} according to $x(t) , y(t) , z(t) $ axes are listed in the Figure \ref{fig:wf}. For the solutions of the system the Runge-Kutta method was employed.

 Obviously the time series show that the functions $x(t)$, $y(t)$, and $z(t)$ are not periodic, indicating that the system is chaotic. 

\begin{figure}[H]
    \centering
    \begin{subfigure}[b]{0.4\textwidth}
        \centering
        \includegraphics[width=\textwidth]{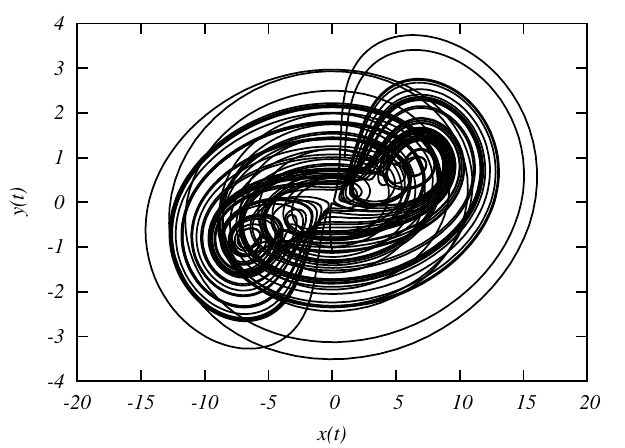}
        \caption{}
        \label{fig:x(t)}
    \end{subfigure}
    \hfill
    \begin{subfigure}[b]{0.4\textwidth}
        \centering
        \includegraphics[width=\textwidth]{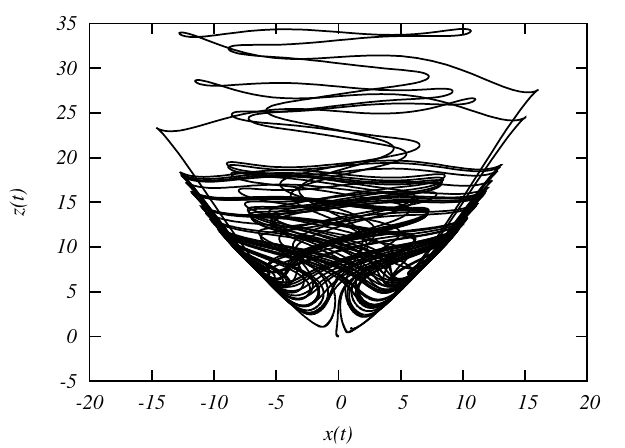}
        \caption{}
        \label{fig:y(t)}
    \end{subfigure}
    \hfill
    \begin{subfigure}[b]{0.4\textwidth}
        \centering
        \includegraphics[width=\textwidth]{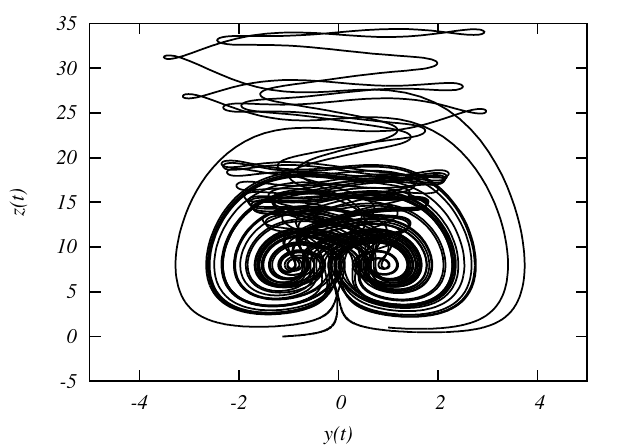}
        \caption{}
        \label{fig:z(t)}
    \end{subfigure}
    \caption{Projection of System \eqref{eq:nlorx}-\eqref{eq:nlorz} on the x-y plane, x-z plane, y-z plane respectively}
    \label{fig:pp}
\end{figure}

The projections of the system \eqref{eq:nlorx}-\eqref{eq:nlorz}, on the various axes is given in the Figure \ref{fig:pp}.

The chaotic attractors are displayed in the Figure \ref{fig:pp}. It appears that the new attractor exhibits an interesting complex chaotic dynamics behavior.

\section{Geometric and Bigeometric Sense of New System}

   As an application of the dynamical systems, the newly formed chaotic system can also be written in the sense of Geometric and the Bigeometric calculus.    

\subsection{Geometric and the Bigeometric Chaotic Systems}

The Geometric Chaotic system corresponding to system \eqref{eq:nlorx}-\eqref{eq:nlorz}, can be written as 

\begin{equation}
\begin{aligned}
\frac{d^*x}{dt}&=e^\frac{\sigma(y z-x)}{x} \\
\frac{d^*y}{dt}&=e^\frac{\rho x-x z}{y} \\
\frac{d^*z}{dt}&=e^\frac{(x y)^2-\beta z}{z} 
\label{eq:newgeo_c}
\end{aligned}
\end{equation}

On the other hand, the Bigeometric Chaotic system corresponding to system \eqref{eq:nlorx}-\eqref{eq:nlorz} will be 

\begin{equation}
\begin{aligned}
\frac{d^\pi x}{dt}&=e^\frac{t(\sigma(y z-x))}{x} \\
\frac{d^\pi y}{dt}&=e^\frac{t(\rho x-x z)}{y} \\
\frac{d^\pi z}{dt}&=e^\frac{t((x y)^2-\beta z)}{z} 
\end{aligned}
\label{eq:newbigeo_c}
\end{equation}

In order to solve the system \eqref{eq:newgeo_c}, one can use the Geometric Runge-Kutta method defined in \cite{RA13}, while the system \eqref{eq:newbigeo_c} can be solved by using the Bigeometric Runge-Kutta method which was defined in \cite{MB}. By choosing the same values for the parameters, such as $\sigma=12$, $\rho=8$ with varying $\beta$ and the initial values as $(0,0,0)$, we can see that the solutions of the Geometric chaotic system \eqref{eq:newgeo_c} and the Bigeometric chaotic system \eqref{eq:newbigeo_c} will be similar to the ones that we get from the solutions of the chaotic system \eqref{eq:nlorx}-\eqref{eq:nlorz}. The following figure shows the chaotic behaviour of the Geometric and the Bigeometric dynamic systems when the parameters are chosen as $\sigma=12$, $\rho=8$ and $\beta=4$.

\begin{figure}[H]

       \centerline{\includegraphics[width=8cm]{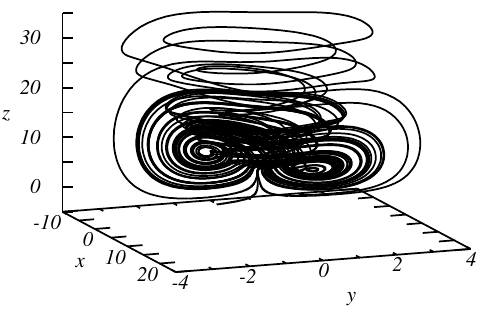}}
        \caption{$\beta=4$}
\end{figure}

\section{Conclusion}

The proposed modified quadratic Lorenz attractor was analysed theoretically and numerically showing that this system is chaotic. Following the standardized analysis method for chaotic systems, we determined first the equilibria points and the eigenvalues of the Jacobian matrix at these equilibria points to get a first indication about the stability of the proposed system. As the eigenvalues are either non-positive real numbers, or complex numbers with positive real parts, we can conclude that this system is not stable at the equilibria points. Moreover, we could identify that the proposed system has a rotational symmetry about the $z$-axis, and shows dissipative behaviour contracting the volume element $V_0$ to $V_0 e^{-16}$. The Lyapunov exponents yield to $l_1 =5.4162$ ,  $l_2 = 2.1912$ and $l_3 = -19.2269$, showing the chaotic nature of the system. The fractional dimension of the system $D_L= 2.3957$ has also been given. Overall, the analysis shows that this is a new chaotic system with two scrolls.

\section{References}

\bibliographystyle{plainnat}

\end{document}